# Rhymes in primes

**Andrei Okounkov**

**Abstract**

While the author is a professional mathematician, he is by no means an expert in the subject area of these notes. The goal of these notes is to share the author's personal excitement about some results of James Maynard with mathematics enthusiasts of all ages, using maximally accessible, yet precise mathematical language. No attempt has been made to present an overview of the current state field, its history, or to place this narrative in any kind of broader scientific or social context. See the references in Section 11 for both professional surveys and popular science accounts that will certainly give the reader a broader and deeper understanding of the material.

## 1. The ancient sieve

It is hard to imagine a more fundamental arithmetic object than the multiplication table

$$\begin{matrix} 1 & 2 & 3 & 4 & 5 & 6 & 7 & \cdots \\ 2 & 4 & 6 & 8 & 10 & 12 & 14 & \cdots \\ 3 & 6 & 9 & 12 & 15 & 18 & 21 & \cdots \\ 4 & 8 & 12 & 16 & 20 & 24 & 28 & \cdots \\ 5 & 10 & 15 & 20 & 25 & 30 & 35 & \cdots \\ 6 & 12 & 18 & 24 & 30 & 36 & 42 & \cdots \\ 7 & 14 & 21 & 28 & 35 & 42 & 49 & \cdots \\ \vdots & \vdots & \vdots & \vdots & \vdots & \vdots & \vdots & \end{matrix} \qquad (1)$$

where the dots indicate that we imagine this table has infinitely many rows and columns. The numbers $n$ that appear in the shaded area are called *composite* numbers. They can be written in the form $n = ab$ where both $a \ne 1$ and $b \ne 1$ are positive integers.

Numbers that are not 1 and not composite are called *prime*. For instance, 2, 3, 5, and 7 are prime, as one sees from (1). Indeed, every composite number $ab$ appears in the multiplication table in the column $a$ and row $b$, which are both *less* than the number $ab$. So, 2, 3, 5, 7 will never appear in the shaded part.

It is a fundamental arithmetic fact that every positive integer $n > 1$ can be factored as a product of primes, and this factorization is unique up to the order of the prime factors. One can compare and contrast factorization into primes with how molecules are built from atoms. One clear difference is that the order of prime factors does not matter, unlike the positions of the atoms in a molecule.

Primes form an infinite sequence which has mesmerized and puzzled mathematicians for millenia. Many mathematicians were first attracted to mathematics by the magic of prime numbers and remained true to their first mathematical love — number theory.

*"It is the fact that primes are so fundamental (being the building blocks of whole numbers), but still so mysterious and poorly understood which makes them so fascinating to me"*, says James Maynard, the hero of these notes. Kannan Soundararajan, the presenter of Maynard's Fields Medal laudatio at ICM 2022, agrees: *"Like many others, I was drawn in by the extreme simplicity of problems involving primes, and the remarkable difficulty of proving anything about them. Twin primes and Goldbach in particular were especially fascinating problems. It's been amazing to witness such spectacular progress as the Green–Tao theorem and bounded gaps between primes over the last twenty years."*

The following method for tabulating the primes goes at least far back as Eratosthenes (276 – 195/194 BC). To remove the composite numbers from the list of all numbers, we can successively cross out or punch trough all numbers from the grey columns in the multiplication table (1), that is, remove all nontrivial multiples of 2, of 3, of 5, et cetera. For instance, the list of natural numbers with 1 and multiples of 2 and 3 removed will look like this:

$$
\begin{array}{cccccccccc}
\bigcirc & 2 & 3 & \cancel{4} & 5 & \cancel{6} & 7 & \cancel{8} & \cancel{9} & \cancel{10} \\
11 & \cancel{12} & 13 & \cancel{14} & \cancel{15} & \cancel{16} & 17 & \cancel{18} & 19 & \cancel{20} \\
\cancel{21} & \cancel{22} & 23 & \cancel{24} & 25 & \cancel{26} & \cancel{27} & \cancel{28} & 29 & \cancel{30} \\
31 & \cancel{32} & \cancel{33} & \cancel{34} & 35 & \cancel{36} & 37 & \cancel{38} & \cancel{39} & \cancel{40} \\
41 & \cancel{42} & 43 & \cancel{44} & \cancel{45} & \cancel{46} & 47 & \cancel{48} & 49 & \cancel{50} \\
\cancel{51} & \cancel{52} & 53 & \cancel{54} & 55 & \cancel{56} & \cancel{57} & \cancel{58} & 59 & \cancel{60} \\
61 & \cancel{62} & \cancel{63} & \cancel{64} & 65 & \cancel{66} & 67 & \cancel{68} & \cancel{69} & \cancel{70} \\
71 & \cancel{72} & 73 & \cancel{74} & \cancel{75} & \cancel{76} & 77 & \cancel{78} & 79 & \cancel{80} \\
\cancel{81} & \cancel{82} & 83 & \cancel{84} & 85 & \cancel{86} & \cancel{87} & \cancel{88} & 89 & \cancel{90} \\
91 & \cancel{92} & \cancel{93} & \cancel{94} & 95 & \cancel{96} & 97 & \cancel{98} & \cancel{99} & 100 \\
\vdots & \vdots & \vdots & \vdots & \vdots & \vdots & \vdots & \vdots & \vdots & \vdots
\end{array}
\tag{2}
$$

where dots indicate that this table has infinitely many rows. The reader may notice there is no need to worry about multiples of 4, 6, or any other composite number.

Once we remove all composite numbers from numbers up to a 100, the result will look like this (the colors will be explained momentarily):

$$
\begin{array}{cccccccccc}
\bigcirc & 2 & 3 & \bigcirc & 5 & \bigcirc & 7 & \bigcirc & \bigcirc & \bigcirc \\
11 & \bigcirc & 13 & \bigcirc & \bigcirc & \bigcirc & 17 & \bigcirc & 19 & \bigcirc \\
\bigcirc & \bigcirc & 23 & \bigcirc & \bigcirc & \bigcirc & \bigcirc & \bigcirc & 29 & \bigcirc \\
31 & \bigcirc & \bigcirc & \bigcirc & \bigcirc & \bigcirc & 37 & \bigcirc & \bigcirc & \bigcirc \\
41 & \bigcirc & 43 & \bigcirc & \bigcirc & \bigcirc & 47 & \bigcirc & \bigcirc & \bigcirc \\
\bigcirc & \bigcirc & 53 & \bigcirc & \bigcirc & \bigcirc & \bigcirc & \bigcirc & 59 & \bigcirc \\
61 & \bigcirc & \bigcirc & \bigcirc & \bigcirc & \bigcirc & 67 & \bigcirc & \bigcirc & \bigcirc \\
71 & \bigcirc & 73 & \bigcirc & \bigcirc & \bigcirc & \bigcirc & \bigcirc & 79 & \bigcirc \\
\bigcirc & \bigcirc & 83 & \bigcirc & \bigcirc & \bigcirc & \bigcirc & \bigcirc & 89 & \bigcirc \\
\bigcirc & \bigcirc & \bigcirc & \bigcirc & \bigcirc & \bigcirc & 97 & \bigcirc & \bigcirc & \bigcirc \\
\vdots & \vdots & \vdots & \vdots & \vdots & \vdots & \vdots & \vdots & \vdots & \vdots
\end{array}
\tag{3}
$$

This table has a lot of holes, just like a sieve. For this reason, the methods that produce an interesting set (e.g. primes) from a less interesting set (e.g. integers) by successively sifting out the unwanted elements are referred to as *sieve* methods.

The primes shown in green are the *twin primes*, that is, primes $p$ such that $p + 2$ or $p - 2$ are also prime[1]. Twin primes are the simplest rhymes in the mysterious poem of primes. While it is very easy to see that there are infinitely many primes[2], the infinitude of twin primes is a very old conjecture, still open today. However, the recent years saw an incredible progress in our understanding of various patterns in primes, recognized, in particular, by the Fields medal, the highest honor in mathematics, awarded in 2022 to James Maynard.

---

[1] Can you prove that $p + 2$ and $p - 2$ cannot both be prime, except for $p = 5$? Questions like this will be clarified when we talk about *admissible* patterns.
[2] Every divisor of the number $n! + 1$, where $n! = 1 \cdot 2 \cdot 3 \cdot \cdots \cdot n$, has to be larger than $n$. Since $n$ is arbitrary, there are infinitely many primes.

3    **Rhymes in primes**

In these notes, we will try to give a very basic introduction to this area of number theory and some of the results of Maynard and his predecessors. A more experienced reader can probably skip many sections of this narrative. All newcomers we wish some patience working through these notes, and very much hope this patience will be rewarded by the sense of awe that this mathematics inspires.

**2. Last digits of primes**

It is very noticeable in (3) that some columns have very few (in fact zero or one) prime numbers in them. Given a number $n$, its column number in (3) is determined by the last digit of $n$ in its decimal notation or, equivalently, by the remainder in the division of $n$ by 10. Mathematicians have a special notation for the remainder, namely

$$89 \bmod 10 = 9.$$

One also says that the *residue* of 89 modulo 10 is 9. More generally, we write

$$a_1 = a_2 \bmod b$$

to mean that $a_1 - a_2$ is divisible by $b$. We say that $a_1$ and $a_2$ are equal mod $b$, or that they are in the same *residue class* modulo $b$.

If $n \bmod 10 = 8$ then $n$ is even and not equal to 2, hence $n$ cannot possibly be prime. Therefore, the 8th column in (3) is empty. Similar reasoning applies to the 2nd, 4th, 5th, 6th, and 10th columns. In due time we will see that prime numbers are approximately evenly distributed among the remaining 4 columns of table (3). Whether the column corresponding to a residue $a$ modulo 10 has many or very few primes is determined by the greatest common divisor $\gcd(a, 10)$. The columns with $\gcd(a, 10) > 1$ contain at most one prime.

The base 10 of the decimal expansion can be replaced by any other base $b > 1$. For instance, $b = 2$ means binary expansions, as exemplified by

$$23 = 10111_{\text{binary}} = 1 \cdot 2^4 + 0 \cdot 2^3 + 1 \cdot 2^2 + 1 \cdot 2^1 + 1 \cdot 2^0. \tag{4}$$

Clearly, for all primes $p \neq 2$ we have $p = 1 \bmod 2$.

Generalizing what we have seen for $b = 10$ and $b = 2$, for any base $b$, primes are approximately evenly distributed among residue classes $a$ modulo $b$ such that $\gcd(a, b) = 1$. The residue classes with $\gcd(a, b) > 1$ contain at most one prime each.

For example, if we replace base $b = 10$ in (3), by $b = 211$, which is a prime number, we will get the following distribution of primes $p \leq 211^2$ (shown by blue or green squares, colors mean the same as in Figure (3)).

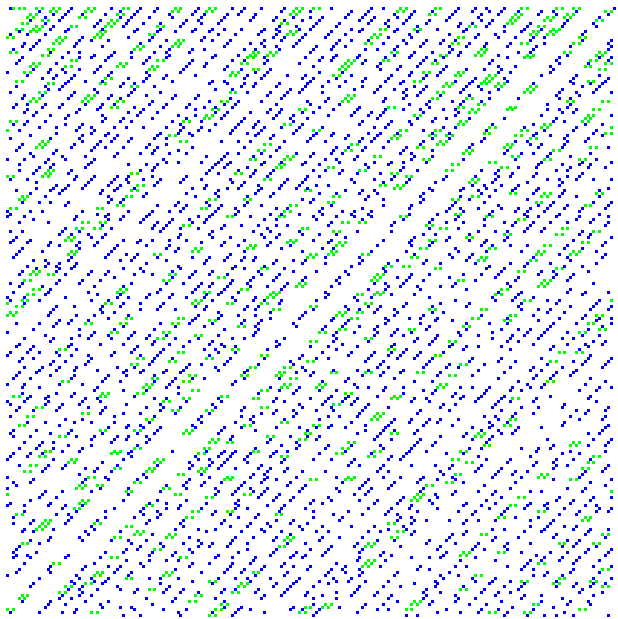

(5)

Primes indeed seem to be roughly evenly distributed among all columns[3], except the very last one, which contains the multiples of 211. Of course, what catches the eye in this picture are the diagonal stripes. We invite the reader to explain them using the equality

$$211i + j = i + j \mod 210$$

and the factorization $210 = 2 \cdot 3 \cdot 5 \cdot 7$.

### 3. The Chinese remainder theorem

One can add and multiply residue classes modulo $b$ in the same way that one can tell the last digit of a sum $n_1 + n_2$ or a product $n_1 n_2$ from the last digits of $n_1$ and $n_2$. Such considerations of are both very basic and very central to number theory. They can be simplified using the Chinese remainder theorem (CRT), which is a result nearly as ancient as the Eratosthenes sieve, appearing in Sunzi Suanjing treatise from the 3rd century CE.

CRT applies to residues modulo $b = b_1 b_2$, where $b_1$ and $b_2$ are *coprime*, meaning that $\gcd(b_1, b_2) = 1$. For example, $10 = 2 \cdot 5$ and $\gcd(2, 5) = 1$. Given a residue $a$ modulo $b$, we can associate to it two numbers

$$a \longrightarrow (a_1, a_2) = (a \mod b_1, a \mod b_2). \tag{6}$$

---

[3] Actually, the number of primes in any given column in (5) varies between 14 and 31, but it all evens out as we go further and further down the list of primes. It is fact of life that it takes a while for primes $p$ to equidistribute mod any fixed prime like $q = 211$. It is a very subtle business to find out how long exactly this while can be, for either some fixed $q$ and or averaged over $q$. This is, in fact, one of the key technical questions in this part of number theory.

For instance, for $b = 10 = 2 \cdot 5$, consider the following table. The rows and the columns of this table are indexed by residues mod 2 and 5, respectively, and we place each residue mod 10 in the corresponding row and column:

|         | 1 | 2 | 3 | 4 | 0 | mod 5 |
|---------|---|---|---|---|---|-------|
| 1 mod 2 | 1 | 7 | 3 | 9 | 5 |       |
| 0 mod 2 | 6 | 2 | 8 | 4 | 0 |       |

$$\qquad(7)$$

We observe the remarkable fact that each residue $a = 0, 1, \ldots, 9$ mod 10 finds a unique place in this table, filling the table completely. In general CRT says that the map (6) gives a one-to-one correspondence

$$\{\text{residues mod } b_1 b_2\} = \{\text{residues mod } b_1\} \times \{\text{residues mod } b_2\} \qquad(8)$$

that preserves arithmetic operations. We invite the reader to prove the CRT and to generalize its statement to the case $b = b_1 b_2 \cdots b_r$.

Let us revisit table (3) from the point of view of CRT. Shading the residue classes that contain $\le 1$ primes, we get

|         | 1 | 2 | 3 | 4 | 0 | mod 5 |
|---------|---|---|---|---|---|-------|
| 1 mod 2 | 1 | 7 | 3 | 9 | 5 |       |
| 0 mod 2 | 6 | 2 | 8 | 4 | 0 |       |

$$\qquad(9)$$

which illustrates two key points:

- $a$ is coprime to 10 if and only if $a$ is coprime to 2 and 5,

- being coprime to 2 and 5 are *independent events*.

Here we think of residue classes $a$ modulo 10 as all equally likely and we call two events $\mathscr{E}_1$ and $\mathscr{E}_2$ independent if

$$\text{Prob}(\mathscr{E}_1 \& \mathscr{E}_2) = \text{Prob}(\mathscr{E}_1) \, \text{Prob}(\mathscr{E}_2) \,.$$

While primes are truly special and not random at all, after centuries of looking into patterns in primes most mathematicians would probably agree that primes behave as if they were completely random, subject to, first, all possible constraints imposed by the considerations of residues and, second, density constraints imposed by the unique factorization of integers into primes. It is therefore very useful to inject, following Cramér, some probabilistic terminology and intuition into our discussion.

### 4. Infinity and limits

There is mystery and challenge in primes because there are infinitely many of them. Any list or plot of primes that we can examine, however long, contains only 0% of all primes, hence always at the best provides a warm-up for the real question. Which is: what happens for all sufficiently large primes?

In mathematics, there is lot of questions for which one is free to discard an arbitrary finite part of some infinite data set. As an example, let's take the concept of a *limit*, which is very important when talking about primes. In the discussion that follows, we will very often have a sequence of real numbers

$$(a_n) = (a_1, a_2, a_3, \ldots),$$

that tends to a limit

$$a = \lim_{n \to \infty} a_n \tag{10}$$

as *n* goes to infinity. Slightly incorrectly, this means that every digit in the decimal expansions of $a_n$'s equals to that of $a$, except for finitely many values of $n$. Any person trained in calculus will be quick to point out some problems with this definition, namely

$$a_n = 10^n \not\to 0,$$

even though every digit of $a_n$ is zero except for one value of $n$, while

$$a_n = 0.\underbrace{999\ldots9}_{n \text{ times}} \to 1.0000\ldots,$$

despite the fact that all displayed digits are different. Readers who are not sure how to fix these issues and feel they could use a more rigorous discussion, can find it in Appendix A.

With the notion of a limit one can define infinite sums and products by

$$\sum_{n=1}^{\infty} a_n = \lim_{N \to \infty} \sum_{n=1}^{N} a_n, \quad \prod_{n=1}^{\infty} a_n = \lim_{N \to \infty} \prod_{n=1}^{N} a_n,$$

when these limits exist. For example, for any number $|x| < 1$, we have

$$x^{\infty} = \lim_{n \to \infty} x^n = 0, \tag{11}$$

and also

$$\sum_{n=0}^{\infty} x^n = \frac{1}{1-x}, \tag{12}$$

which we invite the reader to deduce from (11).

Limits are needed not only for talking about infinite sets, but also as a way to define some very important functions[4]

$$e^x = 1 + x + \frac{x^2}{2} + \frac{x^3}{2 \cdot 3} + \frac{x^4}{2 \cdot 3 \cdot 4} + \cdots = \sum_{n=0}^{\infty} \frac{x^n}{n!}, \tag{13}$$

---

[4] The primary reason the exponential $e^x$ and the natural logarithm $\ln y$ are so important in mathematics is because they solve the simplest *differential* equations, namely $(e^x)' = e^x$ and $(\ln y)' = 1/y$. The reader can check this using the series (13), (15), and the rule $(x^n)' = nx^{n-1}$.

where $e = 2.71828\ldots$ is a famous transcendental number that can be computed by substituting $x = 1$ in the above series. Another important constant that we will meet below is the Euler constant

$$\gamma = \lim_{N\to\infty}\left(\ln N - \sum_{1}^{N}\frac{1}{n}\right) = 0.57721\ldots. \qquad (14)$$

Here and below $\ln y$ denotes the function inverse to (13), which means that by definition

$$\ln e^x = x.$$

It is called the *natural logarithm*, and for arguments in $(0, 2)$ it can be computed using the series

$$\ln(1+y) = y - \frac{y^2}{2} + \frac{y^3}{3} - \frac{y^4}{4} + \cdots = \sum_{n=1}^{\infty}(-1)^{n-1}\frac{y^n}{n}, \quad |y| < 1. \qquad (15)$$

Readers unfamiliar with these functions will discover that the exponential $e^x$ grows very quickly with $x$, making the inverse function $\ln y$ grow very slowly. Notice that the sum in (14), with its minus sign, is the partial sum for $y = -1$ in (15). No wonder it goes to $\ln 0 = -\infty$ as $N$ grows.

While Zeno of Elea (c. 495 – c. 430 BC) made a career out of being confused by the $x = 1/2$ case of (12), we want to stress there are no logical problems whatsoever in thinking about the infinity of primes and about limits. We encourage the reader to embrace these notions as something more true and fundamental than any finite approximations to it.

### 5. The density of primes

If $\mathbb{N} = \{1, 2, \ldots\}$ is the set of natural numbers and $\mathscr{A} \subset \mathbb{N}$ is a subset of it, we define

$$\text{density}(\mathscr{A}) = \lim_{N\to\infty}\frac{|\mathscr{A} \cap \{1, \ldots, N\}|}{N}, \qquad (16)$$

assuming this limit exists. When the limit (16) exists, we will also say that this is the probability that a random natural number is in $\mathscr{A}$.

From table (9) it is clear that

$$\text{density}(\{\text{coprime to } 10\}) = \frac{4}{10} = \frac{1}{2} \times \frac{4}{5}. \qquad (17)$$

Similarly, if $p_1, p_2, \ldots, p_r$ are prime then

$$\text{density}(\{\text{coprime to } p_1 p_2 \cdots p_r\}) = \prod_{i=1}^{r}\left(1 - \frac{1}{p_i}\right). \qquad (18)$$

The equality (18) makes one wonder whether

$$\text{density}(\{\text{primes}\}) \stackrel{?}{=} \prod_{\text{all primes } p}\left(1 - \frac{1}{p}\right). \qquad (19)$$

This is indeed true, but with the clarification that

$$\prod_{\text{all primes } p}\left(1 - \frac{1}{p}\right) \stackrel{!}{=} 0, \qquad (20)$$

as we will see momentarily. Let us look at the reciprocal of the product (18). We have the $x = 1/p$ special case of (12)

$$\frac{1}{1 - \frac{1}{p}} = 1 + \frac{1}{p} + \frac{1}{p^2} + \frac{1}{p^3} + \cdots = \sum_{m \geq 0} \frac{1}{p^m},$$

and multiplying those out for different primes $p_i$, we get

$$\prod_{i=1}^{r} \left(1 - \frac{1}{p_i}\right)^{-1} = \sum_{m_1,\ldots,m_r \geq 0} \frac{1}{p_1^{m_1} p_2^{m_2} \cdots p_r^{m_r}}. \tag{21}$$

If the set $\{p_i\}$ contains all primes that are $\leq N$, then the sum on the right in (21) contains, in particular, the reciprocals of all natural numbers $\leq N$. Therefore, by the existence of the prime factorization, we conclude

$$\prod_{\text{all primes } p \leq N} \left(1 - \frac{1}{p}\right)^{-1} = 1 + \frac{1}{2} + \cdots + \frac{1}{N} + \text{more terms}$$

$$\geq 1 + \frac{1}{2} + \cdots + \frac{1}{N}$$

$$= \ln N + \gamma + o(1), \tag{22}$$

where $\gamma$ is the Euler constant from (14) and $o(1)$ denotes a quantity that goes to 0 as $N \to \infty$. This shows that the rightmost term in

$$0 \leq \text{density}(\{\text{primes}\}) \leq \text{density}(\{\text{coprime to } N!\}) = \prod_{\text{all primes } p \leq N} \left(1 - \frac{1}{p}\right) \tag{23}$$

goes to 0 as $N \to \infty$ and completes the proof of (19).

It is curious to notice that taking logarithms in (22) and using that (15) says that $-\ln(1 - p^{-1}) \approx p^{-1}$ for large $p$, we get

$$\sum_{\text{primes } p} \frac{1}{p} = +\infty. \tag{24}$$

This means that the same computation (22) proves that the density of primes is zero and yet there are sufficiently many primes for the series (24) to diverge, as first noted by Euler.

While we may be disappointed in the fact that the number (19) vanishes, very similar considerations often lead to positive results. For instance, let us consider square-free numbers $n$, that is numbers not divisible by $m^2$ for any $m > 1$. This means

$$n \bmod p^2 \neq 0,$$

for any prime $p$. Referring back to (4), this means that the *two* last digits of $n$ in the expansion base $p$ do not vanish simultaneously. Since this pair of digits is free to take any of the $p^2$ possible values, one can conclude

$$\text{density}(\{\text{squarefree}\}) = \prod_{\text{primes } p} \left(1 - \frac{1}{p^2}\right) = \zeta(2)^{-1} = \frac{6}{\pi^2} \approx 0.6. \tag{25}$$

Here we meet the infinitely famous Riemann $\zeta$-function

$$\zeta(s) = \sum_{n=1}^{\infty} \frac{1}{n^s} = \prod_{\text{primes } p} \left(1 - \frac{1}{p^s}\right)^{-1}, \quad s > 1, \tag{26}$$

and its value $\zeta(2)$ first computed by Euler in 1735. Our earlier computation (20) means that $\zeta(1) = \infty$.

### 6. The prime number theorem

For a set $\mathscr{A}$ of zero density, the numbers (16) go to 0 as $N \to \infty$. A finer measurement of the density is then the rate at which the limit 0 as approached. For prime numbers, the answer is given by the prime number theorem, which says that the density of primes around some large number $N$ is about $1/\ln(N)$.

A mathematically precise way to phrase it uses the function

$$\pi(x) = \text{number of primes } p \text{ such that } p \le x \tag{27}$$

and states that[5]

$$\pi(x) \sim \operatorname{Li}(x) \stackrel{\text{def}}{=} \int_2^x \frac{dy}{\ln y} \sim \frac{x}{\ln(x)}, \tag{28}$$

where $f_1(x) \sim f_2(x)$ means that $\frac{f_1(x)}{f_2(x)} \to 1$ as $x \to \infty$. The reader may find the following data, taken from the Online encyclopedia of integer sequences, convincing:

| $x$ | $\pi(x)$ | $\operatorname{Li}(x)/\pi(x) - 1$ | |
|---|---|---|---|
| 10 | 4 | .25 | |
| $10^2$ | 25 | .16 | |
| $10^3$ | 168 | .054 | |
| $10^4$ | 1229 | .013 | |
| $10^5$ | 9592 | .0039 | |
| $10^6$ | 78498 | .0016 | |
| $10^7$ | 664579 | .00051 | |
| $10^8$ | 5761455 | .00013 | |
| $10^9$ | 50847534 | .000033 | |
| $10^{10}$ | 455052511 | .0000068 | |
| $10^{11}$ | 4118054813 | .0000028 | (29) |
| $10^{12}$ | 37607912018 | .0000010 | |
| $10^{13}$ | 346065536839 | .00000031 | |
| $10^{14}$ | 3204941750802 | .000000098 | |
| $10^{15}$ | 29844570422669 | .000000035 | |
| $10^{16}$ | 279238341033925 | .000000012 | |
| $10^{17}$ | 2623557157654233 | .0000000030 | |
| $10^{18}$ | 24739954287740860 | .00000000089 | |
| $10^{19}$ | 234057667276344607 | .00000000043 | |
| $10^{20}$ | 2220819602560918840 | .00000000010 | |
| $10^{21}$ | 21127269486018731928 | .000000000028 | |
| $10^{22}$ | 201467286689315906290 | .0000000000096 | |

Lest the reader concludes that the last column is always positive, it is known that, in fact, the function $\operatorname{Li}(x) - \pi(x)$ changes sign infinitely many times. Also, while all 3 functions in (28) grow at the same rate, the logarithmic integral Li($x$) gives a much better approximation to $\pi(x)$ than the ratio $\frac{x}{\ln(x)}$.

---

5    A limit procedure is part of the definition of such everyday notions as areas and volumes. The integral of a univariate or multivariate function $f$ is the signed area or volume between the graph of $f$ and the graph of the zero function. It is a continuous limit of summing the values of $f$ over a finer and finer mesh.

The prime number theorem was first shown by Hadamard and de la Vallée Poussin in 1896, so more than 2000 years after Eratosthenes. Certainly, many additional ideas were required, and are still required today to prove (28). Therefore, we will say very little about the proof. The reader interested in a heuristic derivation of the $1/\ln(N)$ density from unique factorization can find it here (requires familiarity with integrals).

To extract the distribution of primes from (93), Hadamard and de la Vallée Poussin had to use some properties of $\zeta(s)$ for complex values of $s$. What happens with $\zeta(s)$ for complex $s$ involves some of deepest problems in all of mathematics, including the infinitely famous Riemann hypothesis (RH), still completely open today. The RH says that all solutions of $\zeta(s) = 0$ are either the so-called trivial zeros $s = -2, -4, -6, \ldots$ or have real part $\Re s = \frac{1}{2}$.

The remarkable $\frac{1}{2}$ from the Riemann Hypothesis can be in fact seen in the table (29) if one notices that the number of 0's in the second column is about half the number of digits of $\pi(x)$, meaning that the difference $\pi(x) - \text{Li}(x)$ of of the order $x^{1/2}$, give or take some logarithmic factors. If there was a zero with $\Re s = c > \frac{1}{2}$, the error $\pi(x) - \text{Li}(x)$ would be at least of size $x^c$, and the argument of Hadamard and de la Vallée Poussin was really about proving that $\Re s < 1$ for all zeros of the $\zeta$-function.

While this is an incredibly interesting topic, the plot of our narrative follows a different path. Asked about the RH, James Maynard says: "*The Riemann Hypothesis suggests that there is a deep hidden structure within the prime numbers. This must occur for a good reason - we just do not know what the reason is yet.*"

### 7. Inclusion–exclusion

Let $\mathscr{A}$ be a set of integers, or even of objects of arbitrary nature. A very, very abstract formulation of a sieve involves some subsets $\mathscr{A}_p \subset \mathscr{A}$, labelled by $p$ in some index set $p \in \mathscr{P}$, which we wish to remove or sift out from the set $\mathscr{A}$. In other words, our goal is to understand the complement $\mathscr{A} \setminus \bigcup_{p \in \mathscr{P}} \mathscr{A}_p$ of all sets $\mathscr{A}_p$ in $\mathscr{A}$.

In its most basic form, the principle of inclusion–exclusion refers to the following elementary observation. Assuming the number of elements $|\mathscr{A}|$ is finite, we have

$$\begin{aligned}
\left|\mathscr{A} \setminus \bigcup_{p \in \mathscr{P}} \mathscr{A}_p\right| = &\ |\mathscr{A}| && \text{count all elements of } \mathscr{A} \\
&- \sum_p |\mathscr{A}_p| && \text{subtract } |\mathscr{A}_p| \text{ for each } p \\
&+ \sum_{p_1 < p_2} |\mathscr{A}_{p_1} \cap \mathscr{A}_{p_2}| && \text{correct for subtracting twice} \\
&- \sum_{p_1 < p_2 < p_3} |\mathscr{A}_{p_1} \cap \mathscr{A}_{p_2} \cap \mathscr{A}_{p_3}| + \ldots && \text{et cetera.}
\end{aligned} \tag{30}$$

For example, referring back to table (9) we may take

$$\mathscr{A} = \{\text{residues modulo } 10\},$$
$$\mathscr{A}_p = \{\text{multiples of } p\}, \qquad p \in \mathscr{P} = \{2, 5\},$$

in which case (30) gives

$$\left|\mathscr{A} \setminus \bigcup_{p=2,5} \mathscr{A}_p\right| = |\{\text{residues coprime to } 10\}| = 10 - 5 - 2 + 1.$$

In other words, subtracting 5 multiples of 2 and 2 multiples of 5, we subtract the zero residue twice, as the shading in table (9) illustrates. Hence we have to put it back.

If the subsets $\mathscr{A}_p \subset \mathscr{A}$ correspond to *independent* events, meaning that

$$\frac{|\mathscr{A}_{p_1} \cap \mathscr{A}_{p_2} \cap \cdots \cap \mathscr{A}_{p_r}|}{|\mathscr{A}|} = \prod_{i=1}^{r} \frac{|\mathscr{A}_{p_i}|}{|\mathscr{A}|}, \tag{31}$$

then formula (30) factors very nicely

$$\frac{|\mathscr{A} \setminus \bigcup_{p \in \mathscr{P}} \mathscr{A}_p|}{|\mathscr{A}|} = \prod_{p \in \mathscr{P}} \left(1 - \frac{|\mathscr{A}_p|}{|\mathscr{A}|}\right), \tag{32}$$

special instances of which we have observed in (17), (18), and (25).

For us, $\mathscr{A}$ will always be some set of integers or residue classes and $\mathscr{A}_d \subset \mathscr{A}$ will denote those divisible by a some number $d$. In this case, all possible intersections in (30) can be described very concretely

$$\mathscr{A}_{p_1} \cap \mathscr{A}_{p_2} \cap \cdots \cap \mathscr{A}_{p_r} = \mathscr{A}_{p_1 p_2 \ldots p_r}, \tag{33}$$

as is illustrated for $\mathscr{P} = \{2, 3, 5\}$ in the following diagram. In (34), we visualize a composite number as a kind of molecule formed by its factors. The primes in $\mathscr{P}$ are assigned three different colors.

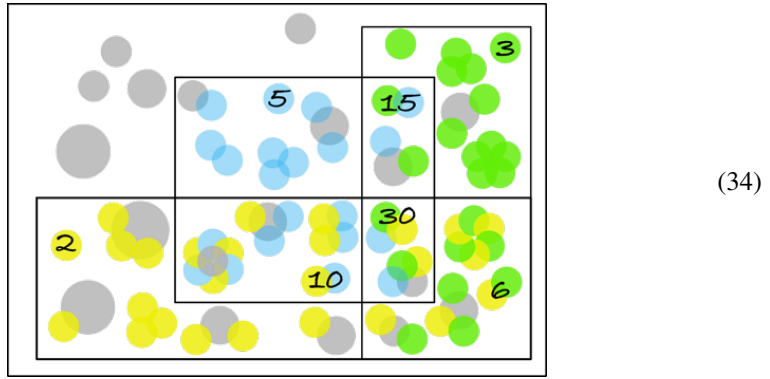

(34)

If (33) is the case, the terms in formula (30) correspond to square-free integers $n$ all prime factors of which belong to $\mathscr{P}$. Thus (30) may be written more compactly

$$\left|\mathscr{A} \setminus \bigcup_{p \in \mathscr{P}} \mathscr{A}_p\right| = \sum_{d=1}^{\infty} \mu_{\mathscr{P}}(d) |\mathscr{A}_d|, \tag{35}$$

using a variant of the Möbius function

$$\mu_{\mathscr{P}}(d) = \begin{cases} (-1)^r, & d \text{ is a product of } r \text{ distict primes in } \mathscr{P}, \\ 0, & \text{otherwise.} \end{cases} \tag{36}$$

A more flexible language for the inclusion-exclusion principle uses the notion of *characteristic functions*. For any subset $S \subset \mathscr{A}$, we define its characteristic function $\delta_S$ by

$$\delta_S(n) = \begin{cases} 1, & n \in S, \\ 0, & n \notin S. \end{cases} \tag{37}$$

Then (35) can be refined to

$$\delta_{\mathscr{A} \setminus \bigcup_{p \in \mathscr{P}} \mathscr{A}_p} = \sum_{d=1}^{\infty} \mu_{\mathscr{P}}(d) \, \delta_{\mathscr{A}_d}. \tag{38}$$

Since

$$|S| = \sum_n \delta_S(n), \tag{39}$$

summing the values in (38) gives (35).

Formulas (35) and (39) require no assumption of indepence like (31). This is very good because (31) is satisfied only approximately in the vast majority of sieve problems. Independence being only approximate is, in fact, a serious problem, to which we will come back below.

Another difficulty one encounters in real number-theoretic applications is that the set $\mathscr{A}$ is typically infinite. For example, we can have $\mathscr{A} = \mathbb{N}$, where $\mathbb{N} = \{1, 2, \dots\}$ is the set of natural numbers. The solution to this problem is to count elements of $n \in \mathscr{A}$ not with weight 1 as in (39), but with some weight $\rho(n)$ such that the count converges. Schematically

$$|\mathscr{A}| = \sum_{n \in \mathscr{A}} 1 \xrightarrow{\text{generalize}} \rho(\mathscr{A}) = \sum_{n \in \mathscr{A}} \rho(n).$$

An example of such weight function is

$$\rho_\zeta(n) = n^{-s}, \quad s > 1, \tag{40}$$

used in the construction of the $\zeta$-function. Multiplicativity of $\rho$

$$\rho(n_1 n_2) = \rho(n_1) \rho(n_2), \tag{41}$$

satisfied by (40) and some other choices of $\rho$, implies an analog of independence (31) for weighted counts. For example, for $\mathscr{A} = \mathbb{N}$, $\mathscr{A}_p = p\mathbb{N}$, and a function $\rho$ satisfying (41), formula (32) transforms into

$$\frac{\sum_{n \text{ coprime to } \mathscr{P}} \rho(n)}{\sum_{n \in \mathbb{N}} \rho(n)} = \prod_{p \in \mathscr{P}} (1 - \rho(p)). \tag{42}$$

We invite the reader to generalize formula (42) for functions $\rho$ satisfying a weaker property

$$\gcd(n_1, n_2) = 1 \implies \rho(n_1 n_2) = \rho(n_1) \rho(n_2). \tag{43}$$

Other than (40), what other interesting functions satisfy (41)? For every $N$, the set

$$(\mathbb{Z}/N\mathbb{Z})^\times = \{\text{residue classes } a \bmod N \text{ such that } \gcd(a,N) = 1\} \qquad (44)$$

is a finite abelian group with respect to multiplication. We take a character of $\chi$ of the group (44) that is, a complex-valued multiplicative function with $\chi(1) = 1$, and extend it by zero to all residues mod $N$. Examples of such functions are

$$\chi_3(n) = \begin{cases} 1, & n = 1 \bmod 3, \\ -1, & n = -1 \bmod 3, \\ 0, & n = 0 \bmod 3, \end{cases} \qquad \chi_5(n) = \begin{cases} i^m, & n = 2^m \bmod 5, \\ 0, & n = 0 \bmod 5, \end{cases} \qquad (45)$$

where the complex number $i = \sqrt{-1} \in \mathbb{C}$ is the imaginary unit. The weight

$$\rho_{N,\chi,s}(n) = \frac{\chi(n \bmod N)}{n^s}, \quad s > 1, \qquad (46)$$

satisfies (41) and the corresponding analog of the $\zeta$-function

$$L(\chi, s) = \sum_{n=1}^\infty \frac{\chi(n \bmod N)}{n^s}, \quad s > 1, \qquad (47)$$

is called the Dirichlet L-function. Its properties are entirely parallel to the $\zeta$-function with one crucial difference. Namely, if $\chi$ is *nontrivial*, that is, takes values other than 0 and 1, then, in contrast to the $\zeta$ function having a singularity at $s = 1$ as in (93), the L-function has a *finite nonzero* value at $s = 1$. This allowed Dirichlet to show that primes are equally distributed among the residue classes (44).

### 8. The first challenge for sieves

As already emphasized above, the main difficulty with sieves is the fact that the independence (31) is only approximate and not exact. Here is an example. Take some large number $x$ and consider the sets

$$\begin{aligned} \mathscr{A} &= \{\text{integers } n \text{ such that } \sqrt{x} < n \le x\}, \\ \mathscr{P} &= \{\text{primes } p \text{ such that } p \le \sqrt{x}\}. \end{aligned} \qquad (48)$$

After sifting out $\mathscr{P}$, we get precisely the primes in the range $(\sqrt{x}, x]$, hence

$$\Big|\mathscr{A} \setminus \bigcup_{p \in \mathscr{P}} \mathscr{A}_p\Big| = \pi(x) - \pi(\sqrt{x}) \sim \frac{x}{\ln x},$$

by the prime number theorem. Let's see if, conversely, we can recover the prime number theorem from the sieve (48).

For fixed $p_1, \ldots, p_r$, the equality (31) is satisfied in the limit $x \to \infty$. However, the error terms present for finite $x$ render the following reasoning *incorrect*. To warn the readers, will use $\stackrel{???}{=}$ to denote an incorrect equality. If we could just apply (32) to the $x \to \infty$ asymptotics, we would get

$$\frac{\pi(x) - \pi(\sqrt{x})}{x - \sqrt{x}} \sim \frac{\pi(x)}{x} \stackrel{???}{\sim} \prod_{\text{primes } p \le \sqrt{x}} \left(1 - \frac{1}{p}\right), \quad x \to \infty. \qquad (49)$$

Having seen products of this general shape before, the reader should not be surprised by the following exact result of F. Mertens

$$\prod_{\text{primes } p \leq \sqrt{x}} \left(1 - \frac{1}{p}\right) \sim \frac{2e^{-\gamma}}{\ln x}, \quad (50)$$

where $\gamma$ is the number from (22) and (93). Since $2e^{-\gamma} \approx 1.123$ this is somewhat close to the right answer and, in particular, gives the correct logarithmic dependence on $x$, but little else can be said in defence of a wrong formula.

This example is meant to illustrate that it is not easy to construct a good sieve, and not to discourage the reader from reading on! See also the references in Section 11, and in particular [7].

## 9. Patterns in primes

So far, we have looked at primes individually, meaning that we studied expressions like

$$\pi(x) = \sum_{\text{primes } p} \delta_{[1,x]}(p), \quad \text{where} \quad \delta_{[1,x]}(y) = \begin{cases} 1, & y \in [1,x], \\ 0, & \text{otherwise}, \end{cases}$$

$$\ln \zeta(s) = -\sum_{\text{primes } p} \ln\left(1 - \frac{1}{p^s}\right),$$

given by summing some natural function $f(p)$ over the set of all primes. To a general science audience, we can say that we have been learning about 1-point correlations in the set of primes.

Recall we expect the primes to be as "random" as the constraints imposed by residues and density allow. To really put these ideas to the test, one should study multi-point correlations, that is, events or patterns that involve pairs, triples, etc. of primes.

To start with a concrete example, what is the probability that $n$ and $n + 1$ are both prime? The answer is clearly 0 because one of these numbers will have to be even, and so $n = 2$ is the only solution. What about $n$ and $n + 2$ being simultaneously prime? Such pairs are called *twin primes* and we saw many such pairs (green) in the Eratosthenes' sieve (3). Similarly, in the plot (5), twin primes are shown in green, all other primes in blue.

Twin primes provide an excellent test of our probabilistic intuition based on density and mod $p$ considerations. From density alone, we should expect that the density of twin primes around $N$ should be about $(\ln N)^{-2}$. However, this needs to be corrected from mod $p$ considerations. Indeed, if $n$ and $n + 2$ were truly independent, the probability of both of them to be coprime to $p$ would be $(1 - 1/p)^2$, while in reality it is $1/2$ for $p = 2$ and $(1 - 2/p)$ for $p > 2$. Whence the following constant in the 1923 conjecture of Hardy and Littlewood

$$\pi_2(x) = |\{p \leq x \text{ such that } p + 2 \text{ is prime}\}| \stackrel{?}{\sim} C_2 \int_2^x \frac{dy}{(\ln y)^2}, \quad x \to \infty, \quad (51)$$

where

$$C_2 = 2 \prod_{\text{primes } p > 2} \frac{1 - \frac{2}{p}}{(1 - \frac{1}{p})^2} = 1.32\ldots. \quad (52)$$

In exactly the same fashion, the probability that $n$ and $n + 2m$ are both coprime to $p$ equals $(1 − 1/p)$ if $p$ divides $2m$ and $(1 − 2/p)$ otherwise. Therefore, for any fixed $m$ one can conjecture that

$$|\{p \le x \text{ such that } p + 2m \text{ is prime}\}| \stackrel{?}{\sim} C_{2m} \int_2^x \frac{dy}{(\ln y)^2}, \quad x \to \infty, \qquad (53)$$

where

$$\frac{C_{2m}}{C_2} = \prod_{p|m,\, p \ne 2} \frac{p-1}{p-2} \ge 1. \qquad (54)$$

From this, it is clear that products of consecutive odd primes like $1155 = 3 \cdot 5 \cdot 7 \cdot 11$ should be particularly likely to occur as distances $p_2 − p_1$ between primes, while powers of two are the least likely values of $p_2 − p_1$. In (55) the function (54) is plotted in the ranges $m \in [1 \ldots 105]$ and $m \in [1 \ldots 1155]$, respectively[6].

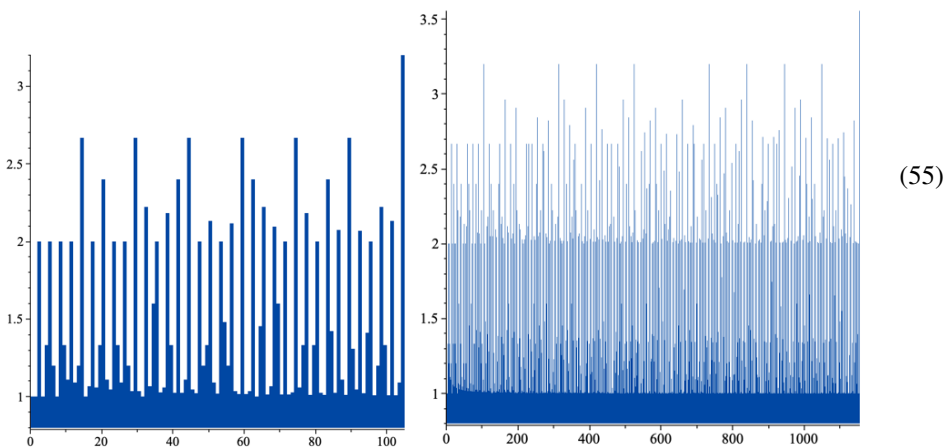

(55)

The conjecture (53) is in excellent agreement with data, especially if one considers the relative frequencies of distances. The following plot (56) compares the function $C_{2m}$ with the actual distribution of the distances among first $10^6$ odd primes:

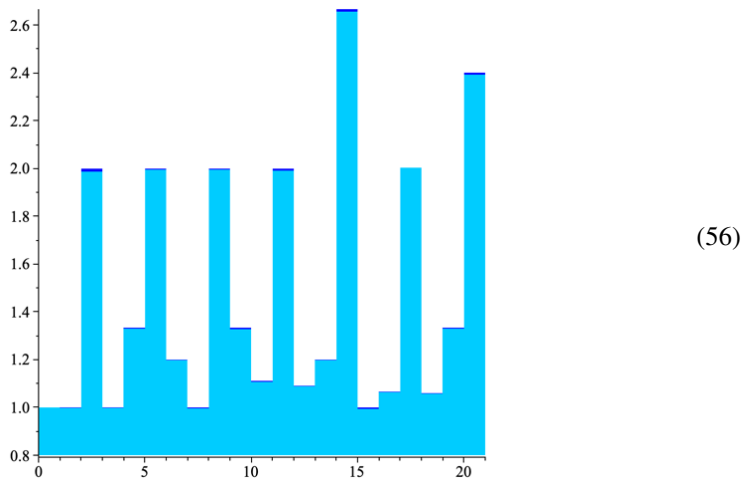

(56)

---

6     The reader may have to adjust the size/resolution of the graph to see the peak at 1155

In (56) we have plotted the relative frequencies, normalized to exactly 1 for $m = 1$. The numerical data is in light blue and the theoretical prediction is in dark blue. The latter overshoots (with the exception of $m = 18$) the former by less than 1%, so it is just barely visible in the plot. Had we gone any deeper in the list of primes, the difference in graphs woud have become undetectable.

We note that the above discussion is for distances between primes, while a *prime gap* of length $2m$ means there are no other primes between $p$ and $p + 2m$. However, since primes become sparser and sparser, finding another prime in an interval of fixed length becomes less and less probable as $p \to \infty$.

The exact same heuristic can be applied to any finite set of jumps

$$J = \{j_1 < j_2 < \cdots < j_l\} \subset \mathbb{N} \tag{57}$$

that we would like to find between primes. We denote by $n + J = \{n + j_1 < \cdots < n + j_l\}$ the shift of $J$ by $n \in \mathbb{N}$ and by $n + J \subset \mathscr{P}$ the event that all numbers $n + j_i$ are prime. In parallel to (53), it is natural to expect that

$$|\{n \le x \text{ such that } n + J \subset \mathscr{P}\}| \stackrel{?}{\sim} C_J \int_2^x \frac{dy}{(\ln y)^{|J|}}, \quad x \to \infty, \tag{58}$$

where

$$C_J = \prod_p \frac{1 - \frac{|J \bmod p|}{p}}{\left(1 - \frac{1}{p}\right)^{|J|}}. \tag{59}$$

Here $|J \bmod p|$ is the number of distinct residue classes mod $p$ in $J$. Since, for fixed $J$, this equals $|J|$ for all sufficiently large $p$, the contribution of all such $p$ to (59) is $1 + O(\frac{1}{p^2})$. Therefore, the product (59) converges.

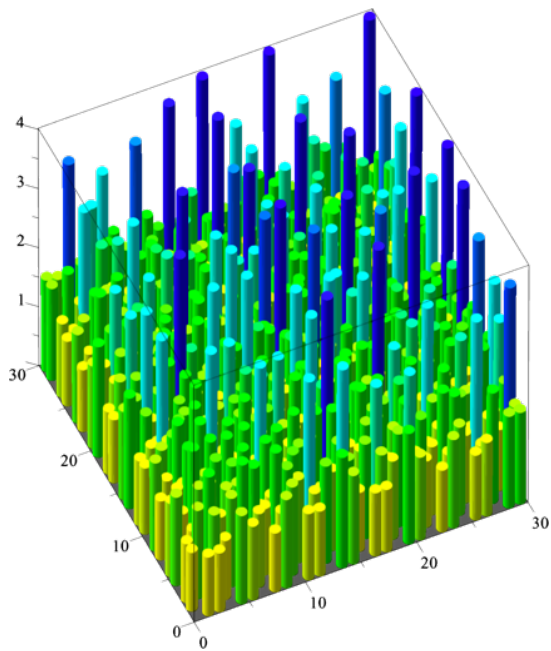

It is clear from (59) that the pattern in primes favor those $J$ that contain a small fraction of residues modulo some prime $p$ and prohibit those $J$ for which $|J \bmod p| = p$. It is also clear from definitions that it suffices to consider the case $j_1 = 0$. The graph of the function $C_{\{0,2i,2(i+j)\}}/C_{\{0,2,6\}}$ is plotted on the left. It vanishes unless $ij(i + j) = 0 \bmod 3$, which explains the missing columns in the plot.

## 10. Closing the gap

Let is call a pattern $J$ as in (57) *admissible* if $C_J \neq 0$, that is, if has a nonzero chance to occur in prime numbers. As a very, very special case of the above heuristic reasoning, one expect that any admissible pattern $J$ will occur as a sequence of prime gaps *infinitely many* times[7]. In particular, one expects the set of twin primes to be infinite. This is known as the *twin prime conjecture*, and it is still open today. However, in constrast to the Riemann Hypothesis, there has been a truly dramatic progress in the recent years on such infinitude questions. This progress has been so dramatic that it inspires us to say that these conjectures are "almost" proven. It is quite incredible to see humans actually reach for the stars.

James Maynard does not quite agree with the narrator here. He says: "*Despite all the recent progress, it seems we are still missing an important idea to prove the Twin Prime Conjecture. But perhaps it is only one big idea.*"

Of course, the actual mathematics involved in proofs compares to what we have discussed so far like a modern airplane compares to a paper airplane. But if the reader tried to think about the issues discussed in Section 8, then she or he may begin to appreciate the amazing creativity and technical mastery required to design sieving arguments leading to the proofs of the breakthrough results below.

It is clear from the prime number theorem that for any constant $c > 1$ there are infinitely many pairs of primes $p_1$ and $p_2$ such that

$$p_1 < p_2 < p_1 + c \ln p_1 . \tag{60}$$

Proving the same statement for some value $c < 1$ is not easy. Many brilliant mathematicians worked on this, finding proofs for smaller and smaller values of $c$, until Goldston, Pintz and Yıldırım have shown that for *any* constant $c > 0$ there are infinitely many pairs of primes satisfying (60).

The new important ideas introduced by Goldston, Pintz and Yıldırım opened the race to replace $c \ln p_1$ in (60) by some fixed constant $B$, that is to prove the infinitude of pairs of primes that are within a fixed finite distance

$$p_1 < p_2 \leq p_1 + B \tag{61}$$

from each other. This race was won in a very dramatic fashion in April 2013 by Yitang Zhang.

Even much more modest results in mathematics today require finding a new way through a real maze of possible ideas, techniques, and logical constructions, and hence moments of extraordinary concentration and clarity of mind. This is not unlike the need to be in a really, really top form for an athlete to set a world record. Research mathematicians (who do have time to do research as part of their job description, in addition to teaching, advising, and other professional duties) cherish these precious moments. Most athletes and mathematicians will surely agree that these special moments tend to be spaced further than $\ln N$ apart once we are past our prime. Zhang's proof is therefore particularly incredible and inspiring, since he had to find his way not just through the mathematical maze, but also through the many turns of his difficult career outside of academia, not giving up despite the big success finally coming to him only at the age of 55. His achievement was widely celebrated by the community, earning him a number of prestigious prizes including the 2013 Ostrowski Prize, the 2014 Cole Prize in Number Theory, and the 2014 Rolf Schock Prize. In the same year 2014, the Cole Prize in Number theory was also awarded to Goldston, Pintz and Yıldırım for their influential work mentioned above.

---

[7] This specific statement is known as the Dickson conjecture, made in 1904.

We hope the reader will turn to [8, 13, 17, 19, 31, 32] to learn more about these developments, and turn to the main hero of these popular notes, the winner of many awards including the 2022 Fields Medal. In the same eventful year 2013, James Maynard realized he can make the sieve a lot more effective, eclipsing Zhang's result in two key dimensions: getting a much stronger result by an easier method.

Speaking about the influences and inspirations that have lead to this result, James Maynard says: "*I was trying to understand the sieve intuition behind the groundbreaking work of Goldston-Pintz-Yıldırım, but in studying this I realised that it might be possible to modify their ideas to go further.*"

It is commonly said that great minds think alike, and the same sometimes happens to the greatest minds, also. In the suspenseful race to close the prime gap, Terry Tao arrived at the same results independently at the same time as James Maynard. "*I was a bit shocked when I first heard the news, but fortunately Tao was very generous and understanding. Simultaneous discovery happens more often than you'd imagine*!", says James Maynard.

To explain Maynard's and Tao's main result on small gaps in primes, it is important to make a certain change of perspective. In Section 9, we were interested in the event when *all* numbers

$$n + J = (n + j_1, n + j_2, \ldots, n + j_l) \tag{62}$$

are prime. But if one asks for less one can prove more! Let's instead fix some $m < l$ and ask that at least $m$ of the numbers (62) are prime for infinitely many values of $n$. We will not know which ones among (62) are prime, but we will know, for instance, that there are infinitely many primes within distance $j_l - j_1$ from each other.

The following is a special case of the spectacular main result of [20], which Kannan Soundararajan compares with "*sun amidst the stars*" in his Fields Medal laudatio.

**Theorem 1.** *For any m, for all sufficiently long admissible patterns J, at least m of the numbers* (62) *are prime for infinitely many n.*

In fact, for any given *m*, the required size of *J* in Theorem 1 can be made explicit. For $m = 2$, $|J| = 50$ suffices, and the following set being admissible

$$\begin{aligned} J = \{&0, 4, 6, 16, 30, 34, 36, 46, 48, 58, 60, 64, 70, 78, 84, 88, 90, 94, 100, 106, \\ &108, 114, 118, 126, 130, 136, 144, 148, 150, 156, 160, 168, 174, 178, 184, \\ &190, 196, 198, 204, 210, 214, 216, 220, 226, 228, 234, 238, 240, 244, 246\}, \end{aligned} \tag{63}$$

shows there are infinitely many primes at most 246 apart.

For $m = 3$, $|J| = 35410$ suffices, and one can take[8], for instance, the first 35410 primes larger than 35410

$$J = \{35419, 35423, \ldots, 469411, 469397\}\,.$$

Therefore, there are infinitely many triples of primes within 433992 of each other. In general, the best estimate for required length of $J$ currently stands at $ce^{3.815m}$, see [1].

The more general result proven in [20] guarantees there are at least $m$ primes among the numbers $a_1 n + j_1, \ldots, a_l n + j_l$ provided these are distinct and admissible. This stronger version of Theorem 1 leads to many further interesting conclusions about patterns in primes. For example, one can deduce that there are arbitrarily large sets of primes where any pair in the set differs in only 2 decimal places! Indeed, if we take

$$a_i = l!\, 10^{l+2}\,, \quad j_i = 10^{i+1} + 1\,, \tag{64}$$

then all digits of $a_i n + j_i$, $i = 1, \ldots, l$ are the same, except the position of the 1 in the $(i+1)$st decimal place, which is changing its position within the string of $l$ zeros.

I hope the readers share the narrator's sense of awe at this absolutely amazing mathematics and join me in warmest congratulations on it being recognized by the Fields Medal. I also hope the readers got the sense that today's mathematics is not just extraordinarily powerful, but also concrete, understandable, and fun, once one finds the right idea and the right point of view. While finding that right point of view is not at all easy, my biggest hope is to have inspired my youngest readers to believe that mathematics can be beautiful and rewarding, both as a subject and as a profession. Maybe this is also a good place for me to thank James Maynard and Kannan Soundararajan for this special opportunity to be introduced to their wonderful subject.

### 11. Further reading

The *Quanta Magazine* has published several popular accounts of these and related developments, see [11, 13–15, 19].

Among surveys written by top experts in the field, one should mention [5, 8, 17, 27], including expositions by James Maynard himself [21–23].

Among textbooks of different level, the reader will surely find something which suits her or his level and style among [3, 9, 16, 28, 29] or the more advanced [4, 12]. There is even a graphic detective novel [10]!

I hope the reader has a lot of fun studying these sources as well as the original articles [6, 20, 25, 26, 32].

---

[8] As an exercise, the reader may check than any $l$-tuples of primes larger than $l$ is admissible

## 12. A glimpse into the argument

To help the reader make transition to further popular and research reading, we will indicate some initial logical steps in the argument leading to the proof of Theorem 1. There is a certain distance that we can fly even on our paper airplane.

### 12.1. Being prime on average

We need to prove that at least $m$ of the numbers (62) are prime for infinitely many $n$. Suffices to show that for any given integer $N$ this is true for some $n \geq N$. Let $\mathscr{P}$ denote the set of all primes. Instead of trying to find a specific $n$ for which the intersection $\{n + J\} \cap \mathscr{P}$ has at least $m$ elements, we can ask about the average size of the intersection $\left|\{n + J\} \cap \mathscr{P}\right|$ with respect to some density $\rho(n) \geq 0$ on $[N, \ldots, 2N]$. This density $\rho$ is something we are *bringing* into the argument, not something given to us in advance.

Clearly,

$$\text{average}\left(\left|\{n+J\} \cap \mathscr{P}\right|\right) = \frac{\sum \rho(n)\left|\{n+J\} \cap \mathscr{P}\right|}{\sum \rho(n)} \leq \max\left(\left|\{n+J\} \cap \mathscr{P}\right|\right), \quad (65)$$

and so if we can bound the average in (65) below by $m$ then we win. Now, since the numbers $j_k \in J$ are all distinct, we have

$$\frac{1}{\sum \rho(n)} \sum_{n=N}^{2N} \rho(n)\left|\{n+J\} \cap \mathscr{P}\right| = \sum_{k=1}^{l} \frac{\sum_{n+j_k \text{ is prime}} \rho(n)}{\sum_{N \leq n \leq 2N} \rho(n)}. \quad (66)$$

Hence, our strategy is to invent a function $\rho(n)$ for which each of the $l$ ratios in the right-hand side of (66) can be shown to be large.

### 12.2. Looking for $\rho$, part I

A naive strategy would be to take

$$\rho_0(n) = \begin{cases} 1, & n + J \subset \mathscr{P}, \\ 0, & \text{otherwise}. \end{cases} \quad (67)$$

This makes the numerator and denominator in (66) equal, and so naively each fraction equals 1. What this overlooks is that $\frac{0}{0}$ is no good in (66), and that our original goal is precisely equivalent to showing that $\rho_0$ takes some nonzero values.

This underscores the point that we haven't really advanced on the problem yet, just put in a slightly more flexible framework by introducing the density $\rho$. Those who can design a good $\rho$ are the great masters of the sieve.

Functions that only take values 0 or 1 are called *characteristic* functions as we recall from (37). These are also the functions that are equal to their own square. From the definitions,

$$\rho_0(n) = \delta_{[N,\ldots,2N]}(n) \prod_{k=1}^{l} \delta_{\mathscr{P}}(n + j_k). \quad (68)$$

The next natural idea is to find a working replacement $\widetilde{\delta}$ for $\delta_{\mathscr{P}}$ and get $\rho$ by multiplying them together.

Plots of the function $\delta_\mathscr{P}$ look like barcodes, and here is an example

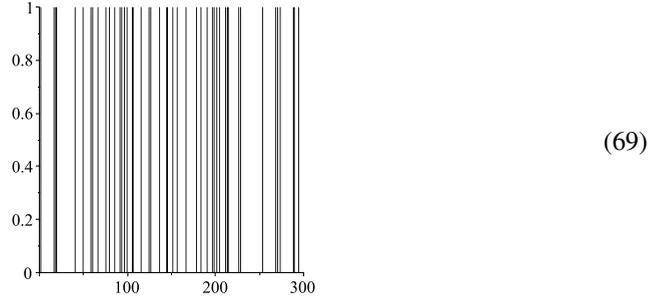
(69)

in which $n$ takes odd values from $10^6 + 1$ to $10^6 + 599$. In principle, (38) gives a formula for $\delta_\mathscr{P}$, and we can approach the goal of finding a replacement $\delta_\mathscr{P}$ by tinkering with the formula (38). For instance, we just truncate summation over $d$ to some maximal value $D$. That is, we define

$$\widetilde{\delta}_0(n) = \left( \sum_{d|n,\, d\leq D} \mu(d) \right)^2, \qquad (70)$$

where we square the sum to make the result nonnegative. Since this equals 1 if $n$ has no nontrivial divisors $d \leq D$, it is natural to compare this function to the characteristic function $\delta_{\leq D}$ of numbers without prime factors $p \leq D$.

It is easy to plot the function $\widetilde{\delta}_0 - \delta_{\leq D}$ and the result

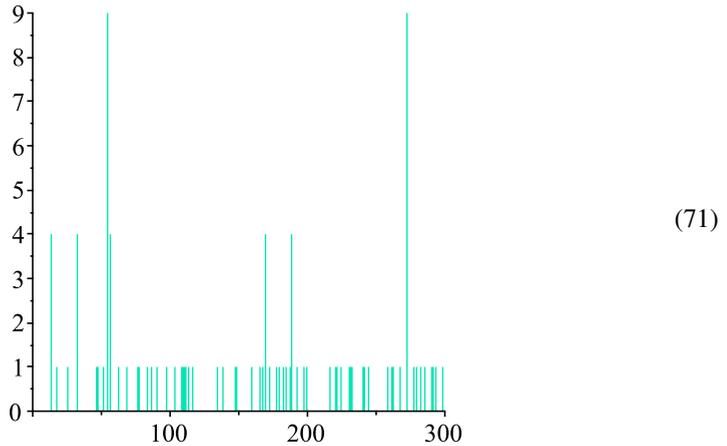
(71)

for $D = 100$ is not really satisfying. The two peaks in the graph correspond to the numbers

$$1000109 = 11 \cdot 23 \cdot 59 \cdot 67, \quad 1000545 = 3 \cdot 5 \cdot 7 \cdot 13 \cdot 733,$$

and, in general, the function (70) becomes large not because $n$ is prime, but because there is a significant disbalance between its divisors $d \leq D$ with different parity of the number of prime factors. In other words, $\widetilde{\delta}_0(n)$ is much more sensitive to the artificial cutoff introduced by us at $d \leq D$ than to what we set out to measure in the first place.

To get rid of this effect, it makes sense to replace the hard cutoff at $d \leq D$ by a more gentle one, through some weight function of $d$ that gives 1 for prime numbers and vanishes at $d = D$. Let us try

$$\widetilde{\delta}_k(n) = \frac{1}{(\ln D)^{2k}} \left( \sum_{d|n,\, d \leq D} \mu(d) \left( \ln \frac{D}{d} \right)^k \right)^2, \tag{72}$$

and this works much, much better for $k \geq 1$. For $D = 100$, the function $\widetilde{\delta}_1 - \delta_{\leq D}$ looks like this:

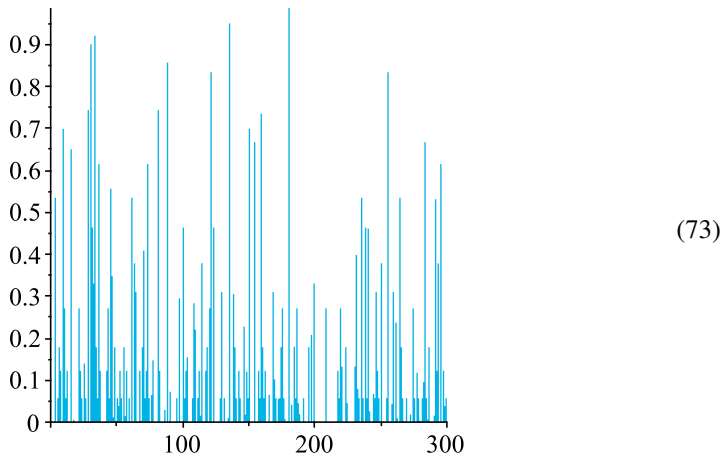

(73)

Not only it takes values in $[0, 1)$ in this plot, it also peaks at numbers with prime factors $p$ of size close to $D$. Since the weight $\ln \frac{D}{p}$ gets small for such $p$, we certainly expect such numbers to contribute on par with the prime numbers.

### 12.3. Looking for $\rho$, part II

Functions (72) played an important role in the work of Goldston, Pintz and Yıldırım. However magical, by themselves they are not enough to get to the Maynard-Tao theorem. If we just multiply them as in (68), then we loose the crucial synergy between different elements of the list $J$. Recall that the logic of Theorem 1 is such that the longer the list $J$ gets, the easier it is to find many prime numbers in it. For this, there should be some nontrivial interaction between different $j_k$.

One key new ingredient in the Maynard-Tao method is to consider functions of the form

$$\rho(n) = \delta_{[N,\ldots,2N]} \left( \sum_{\substack{d_1|n+j_1,\ldots,d_l|n+j_l, \\ d_1 d_2 \cdots d_l \leq D}} \mu(d_1 d_2 \cdots d_l) F\left( \frac{\ln d_1}{\ln D}, \ldots, \frac{\ln d_l}{\ln D} \right) \right)^2, \tag{74}$$

where $F$ is a multivariate function to be specified later. As before, we want $F$ to be small if the arguments sum to 1 (meaning that $d_1 d_2 \cdots d_l = D$) to soften the effect of the summation cutoff introduced in (74).

By allowing $F$ to depend on each divisor $d_i$, the Maynard-Tao method activates a very powerful principle of measure concentration in high-dimensional geometry. At the risk of being repetitive, one may note that there is really a lot of space in a space of a large dimension $N$. There is so much space that no probability distribution can cover all of it evenly as $N \to \infty$, and one could put this vague principle in a mathematically precise form, see for instance [18].

To make a negative statement positive, one can say that any high-dimensional probability density has to concentrate on some small portion of the whole space. For example, a probability measure $\nu$ on the line $\mathbb{R}$ is another name for a random variable $x$, and a product measure $\nu^{\otimes N} = \nu \times \cdots \times \nu$ on $\mathbb{R}^N$ is another name for a sequence of independent, identically distributed (i.i.d.) random variables $x_1, \ldots, x_N$. We know from basic probability theory that, with minimal assumptions about $\nu$, the average $\frac{1}{N} \sum x_i$, and many other functions of i.i.d. random variables $x_1, \ldots, x_N$, will sharply peak, or concentrate, around their expected value as $N \to \infty$.

A reader not familiar with these notions, may experiment by working out the example in which $\nu$ is the uniform density on $[0, 1]$ and $\nu^{\otimes N}$ is a uniform density on an $N$-dimensional cube $[0, 1]^N$. Taking the sum $\sum x_i$ means projecting the cube onto the $(1, 1, \ldots, 1)$ axis, and the reader may enjoy actually plotting these densities for different values of $N$. It is also fun to compute the projection of a uniform measure on a high-dimensional sphere onto any axis.

It is by harnessing these concentration of measure phenomena that the density (74) can significantly improve upon (72).

### 12.4. Primes in arithmetic progressions, on average

Now let's plug the formula (74) into the numerator in (66), expand out the square, and do summation over the variable $n$ first. We get a sum of the form

$$\sum_{n+j_k \text{ is prime}} \rho(n) = \sum_{\vec{d},\vec{d}'} \mu\mu FF \sum_{\text{certain } n} 1 \qquad (75)$$

where the outer sum is over two sets of integers

$$\vec{d} = (d_1, \ldots, d_l) \quad \text{and} \quad \vec{d}' = (d'_1, \ldots, d'_l),$$

there is a weight of the form

$$\mu\mu FF = \mu(\Pi d_i)\mu(\Pi d'_i) F\left(\tfrac{\ln \vec{d}}{\ln D}\right) F\left(\tfrac{\ln \vec{d}'}{\ln D}\right)$$

and the inner sum runs over $n$ such that

$$n + j_i \ = 0 \bmod \mathrm{lcm}(d_i, d'_j), \quad i = 1 \ldots, l, \qquad (76)$$
$$n + j_k \ \text{ is prime}, \qquad (77)$$

where $\mathrm{lcm}(d_i, d'_j)$ denotes the least common multiple.

It is clear from this that we must have $d_k = d'_k = 1$. Since the remaining congruence conditions can be put into a single congruence condition using the Chinese Remainder Theorem, the sum over $n$ thus counts primes in an arithmetic progression.

Time and time again in these notes we have stressed the technical importance of being able to accurately count primes in arithmetic progression in analytic number theory, also stressing that this may be very delicate if the progression is not much longer than its common difference.

The counting function (27) may be refined to count primes in a given residue class modulo $b$

$$\pi(x, b, a) = \text{number of primes } p \text{ such that } p \leq x \text{ and } p = a \bmod b. \tag{78}$$

The Dirichlet theorem mentioned in Section 7 says that

$$\frac{\pi(x, b, a)}{\pi(x)} \to \begin{cases} \phi(b)^{-1}, & \gcd(a, b) = 1, \\ 0, & \text{otherwise}, \end{cases} \tag{79}$$

as $x \to \infty$, where $\phi(b)$ is the number of residue classes coprime to $b$. For fixed $x$, however, the function

$$(b, a) \mapsto \phi(b) \frac{\pi(x, b, a)}{\pi(x)} - 1 \tag{80}$$

behaves in a very irregular manner. This is illustrated in the following plot for $a < b \leq 100$

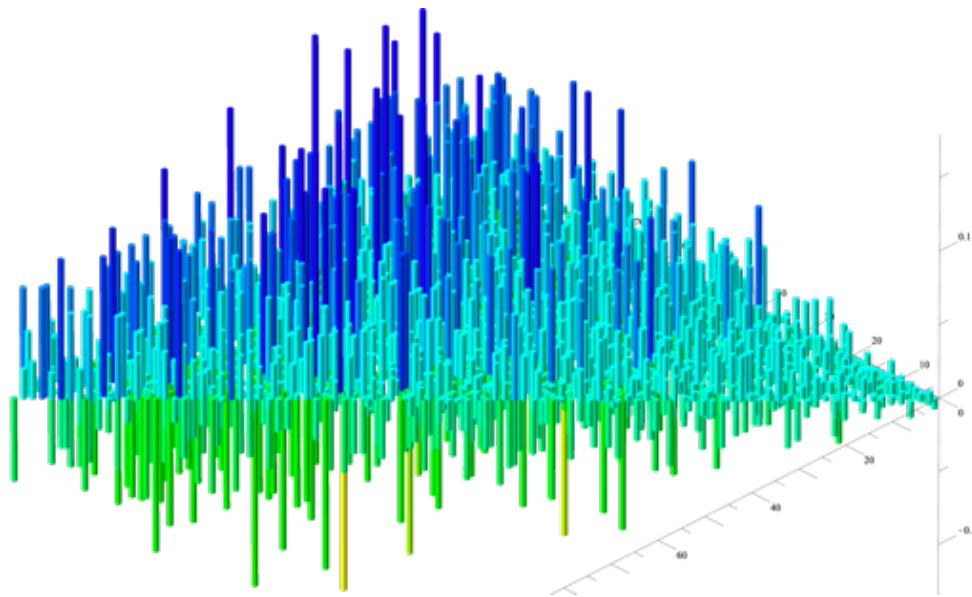

and the first 5000 primes, which means $x^{1/2} \approx 220$.

Very fortunately, in (75), we don't have to face the full complexity of this function. Since there is an outside summation over $\vec{d}$ and $\vec{d}'$, we only need to know its *average* over $b$.

Recall that the Riemann hypothesis implies error of size about $x^{1/2}$ in the prime number theorem. The conjectural extension of the Riemann hypothesis to Dirichlet L-functions (47) would give a similar error bound for $\pi(x, b, a)$. If one sums these errors for $b < x^{1/2}$, one thus expects to get something of order $x$. Remarkably, a slight weakening of this statement, known as the Bombieri-Vinogradov theorem *has* been proven [2, 30]. In other words, the Riemann hypothesis for L-functions is a complete mystery, but its main consequence for the distributions of primes in arithmetic progression can be rigorously proven *on average*. The actual estimate one needs here has the form

$$\sum_{b < x^{1/2-\varepsilon}} \max_{\gcd(a,b)=1} \left| \pi(x,b,a) - \frac{\pi(x)}{\phi(b)} \right| \le C(A,\varepsilon) \frac{x}{(\ln x)^A}, \quad (81)$$

which holds for any $A > 0$ and $\varepsilon > 0$ with some positive constant $C(A,\varepsilon)$ that depends on $A$ and $\varepsilon$. In our example, the maxima over $a$ in (81) and their running average over $b$ can be seen in the following plot

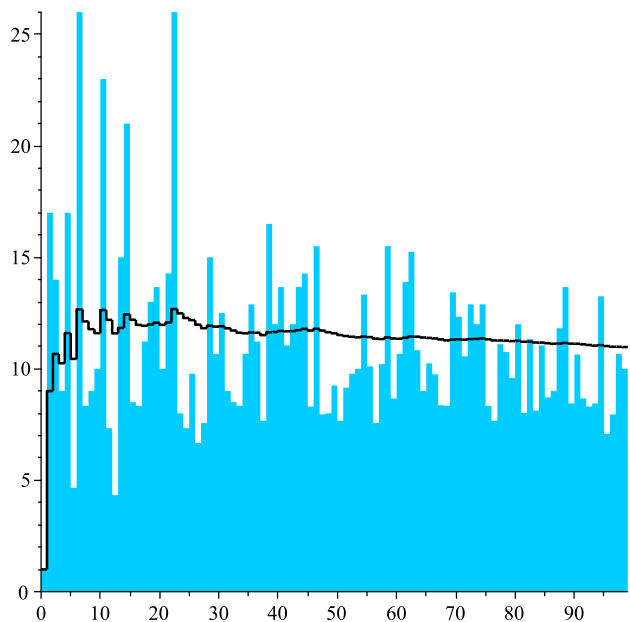

Averaging really does make the behavior a lot more regular and, hence, manageable.

We have discussed some of the key ingredient that go into the proof of the amazing result of Maynard and Tao. Perhaps, this discussion has given the reader the motivation and confidence to open more advanced literature written by the experts in the field, including the papers listed in Section 11. In any case, we hope to have communicated to the reader our own sense of awe at the beauty of mathematics.

## A. Limits

Limits are defined not just for numerical sequences $(a_1, a_2, \dots)$ but for objects of arbitrary nature for which there is a notion of *neighborhoods*. Namely, $a$ is the limit of the above sequence, if *every* neighborhood of $a$ contains all elements $a_n$ except maybe finitely many. The reader may find it useful to picture this as follows:

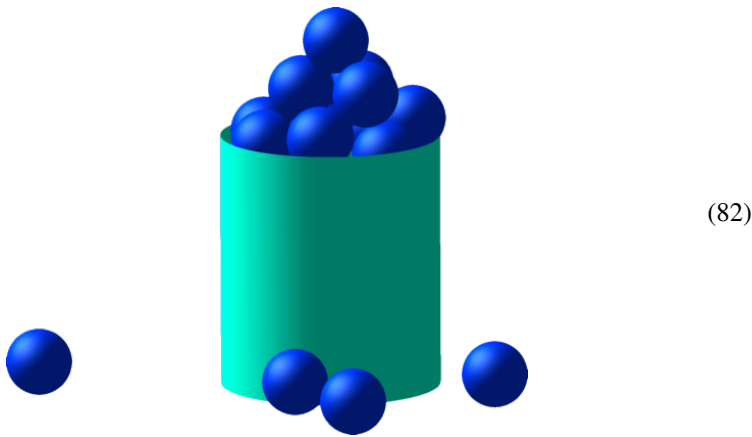

(82)

where the bin represents a neighborhood of $a$ and spheres represent the elements $a_n$. Of course, since the sequence is infinite, any neighborhood of the limit point contains not just many, but infinitely many of the $a_n$'s.

For real numbers, or any other set with the notion of distance, we may take the open balls of arbitrary positive radius $r > 0$

$$B(a, r) = \{\text{all } x \text{ such that distance}(x, a) < r\}$$

as standard neighborhoods. The reader may check her or his understanding of the definition by proving (11) and (12), constructing a sequence or real numbers that does not have a limit, and proving that the limit of a sequence of real numbers is unique when it exists.

The slight issue with defining the limits digit by digit is that the set of all real numbers whose decimal expansion is fixed up to a certain point is a half-open interval, for instance

$$\{\text{all } x \text{ such that } x = 2.71 \dots\} = [2.71, 2.72) \,.$$

To define limits for real numbers correctly, one one should take *open* intervals, that is, those without both endpoints as neighborhoods. Back to the main text.

### B. Mellin transform and the density of primes

Consider a simplified model, in which we forget about integrality and talk about real numbers $x > 1$. Let $\rho_1(x)$ be a certain density function on $[1, \infty)$. It will model the density of prime numbers. What should then correspond to the density $\rho_r(y)$ of the numbers $y$ that have exactly $r$ prime factors?

We have, by definition, $y = x_1 x_2 \dots x_r$, where $x_i$ are distributed in the set

$$\{1 \leq x_1 \leq x_2 \leq \dots \leq x_r\}$$

with density $\rho_1(x_1) \cdots \rho_1(x_r)$. Thus for *any* function $f(y)$ we have

$$\int f(y) \rho_r(y) \, dy = \int_{1 \leq x_1 \leq x_2 \leq \dots \leq x_r} f(x_1 \cdots x_r) \prod \rho_1(x_i) \, dx_i \,. \tag{83}$$

Which functions $f(y)$ should we consider?

In mathematics, the success often depends on choosing the right point of view. If one has the right point of view, then one is able to see clearly where one is going.

A very nice choice here is to take $f(y) = y^{-s}$, where $s > 1$ is parameter. This is called *Mellin transform*, and it is a transform because it takes a function $\rho_r(y)$ of one variable $y$ to another functions $\rho_r^{\text{Mellin}}(s)$, of the parameter $s$. Thus one trades a function of one variable $\rho_r(y)$ for another function of one variable $\rho_r^{\text{Mellin}}(s)$, which seems like a fair exchange. In fact, one can reconstruct $\rho_r(y)$ from $\rho_r^{\text{Mellin}}(s)$, so no information is lost.

The Mellin transform is a close relative of the Fourier transform[9] and what makes the following computation work is the basic identity

$$(x_1 x_2)^s = x_1^s x_2^s .$$

Because of this, the function $f(x_1 \cdots x_r)$ in (83) factors as $f(x_1) \cdots f(x_r)$ and we can eventually reduce an $r$-fold integral in (83) to a product of $r$ integrals.

We compute

$$\rho_r^{\text{Mellin}}(s) \stackrel{\text{def}}{=} \int_1^\infty y^{-s} \rho_r(y)\, dy \tag{84}$$

$$= \int_{1 \le x_1 \le x_2 \le \cdots \le x_r} (x_1 \cdots x_r)^{-s} \prod \rho_1(x_i)\, dx_i \tag{85}$$

$$= \frac{1}{r!} \int_{[1,\infty)^r} \prod x_i^{-s} \rho_1(x_i)\, dx_i \tag{86}$$

$$= \frac{1}{r!} \rho_1^{\text{Mellin}}(s)^r , \tag{87}$$

where in going from (85) to (86) we used the fact that

$$[1,\infty)^r = \bigcup_{\substack{\text{permutations} \\ w: \{1,\ldots,r\} \to \{1,\ldots,r\}}} \{1 \le x_{w(1)} \le x_{w(2)} \le \cdots \le x_{w(r)}\} \tag{88}$$

and that the integration over any of the $r!$ sets in the right-hand side of (88) gives the same result as (85).

If $\rho_\bullet$ is the density of numbers $y$ having an arbitrary number of factors $r$, including the case when $r = 0$ and $y = 1$, then summing (87) over $r = 0, 1, 2, \ldots$ gives

$$\rho_\bullet^{\text{Mellin}}(s) = \exp\left(\rho_1^{\text{Mellin}}(s)\right) , \tag{89}$$

where $\exp(x)$ is another notation for the function $e^x$ from (13). The appearance of the exponential function here is typical in many inclusion-exclusion situations.

To model unique factorization we want to take $\rho_\bullet = 1$ on $[1,\infty)$ which means

$$\rho_\bullet^{\text{Mellin}}(s) = \int_1^\infty x^{-s} dx = \frac{1}{s-1}, \quad s > 1 . \tag{90}$$

Thus, we expect

$$\int_1^\infty x^{-s} \rho_1(x)\, dx \stackrel{?}{=} \ln \frac{1}{s-1}, \quad s > 1 , \tag{91}$$

---

 **9**   Some reader may find the explanation of Fourier transform in [24] usable.

which is both good and bad news for the following reasons.

On the one hand, $\ln \frac{1}{s-1}$ is not a Mellin transform of any density $\rho_1$ on $[1, \infty)$ simply because it does not have a limit as $s \to +\infty$. The $s \to +\infty$ limit in (91) probes $\rho_1(x)$ for $x$ very close to 1 because $x^{-s}$ becomes very small on the whole interval $(1 + \delta, \infty)$ as $s \to \infty$, for any fixed $\delta > 0$. In particular, the Mellin transform of a bounded density function $\rho_1(x)$ on $[1, \infty)$ has to go to zero as $s \to +\infty$.

This means that we cannot accurately model prime numbers with real numbers and continuous densities. Of course, it was certainly silly to be asking for the density of small primes to begin with. However, our interest is precisely the opposite, as we want to know the behavior of $\rho_1(x)$ for large $x$. This region is probed by $s \to 1$ limit of the Mellin tranform. In fact

$$f(x) = f_0 + O(x^{-c}) \Rightarrow \int_1^\infty f(x) x^{-s} dx = \frac{f_0}{s-1} + \ldots, \qquad (92)$$

where $O(x^{-c})$ means that $\left|\frac{f(x) - f_0}{x^{-c}}\right|$ remains bounded as $x \to \infty$, the double arrow $\Rightarrow$ denotes implication, and dots stand for a function which is analytic for $s > 1 - c$. (And also analytic for complex values of $s$ such that $\Re s > 1 - c$.) In the $s \to 1$ limit, we may write

$$\int_1^\infty x^{-s} \rho_1(x) \ln(x) dx = -\frac{d}{ds} \int_1^\infty x^{-s} \rho_1(x) dx \sim -\frac{d}{ds} \ln \frac{1}{s-1} = \frac{1}{s-1}.$$

which strongly suggests $\rho_1(x) \sim 1/\ln(x)$ for $x \to \infty$.

In place of continuous approximations, the proof of Hadamard and de la Vallée Poussin uses properties of the $\zeta$-function (26), which, in the spirit of (84), can be interpreted as the averaged value of $n^{-s}$ with respect to the measure that gives every positive integer $n$ weight 1. The equality between the sum and the product in (26) is the correct discrete version of the relation (89). It looks different because in the discrete situation we need to account for the nonzero chance of having two equal prime factors, the possibility of which was ignored in going from (85) to (86). The exact analog of (90) is the the following description

$$\zeta(s) = \frac{1}{s-1} + \gamma + o(1), \quad s \to 1, \qquad (93)$$

of the $s \to 1$ behavior of the $\zeta$-function, where $\gamma$ is the constant from (14) and (22). Back to the main text.

**Andrei Okounkov**

Andrei Okounkov, Department of Mathematics, University of California, Berkeley, 970 Evans Hall Berkeley, CA 94720–3840, okounkov@math.columbia.edu